\documentclass[12pt,a4paper]{article}
\usepackage{amsmath}
\usepackage{amsfonts}
\usepackage{amssymb}
\usepackage[T1]{fontenc}
\usepackage[ansinew]{inputenc}
\usepackage{latexsym}

\newtheorem{theorem}{Theorem}[section]
\newtheorem{lemma}[theorem]{Lemma}
\newtheorem{remark}[theorem]{Remark}
\newtheorem{remarks}[theorem]{Remarks}
\numberwithin{equation}{section}
\newtheorem{corollary}[theorem]{Corollary}

\def\N{\mathbb{N}}
\def\R{\mathbb{R}}
\def\Q{\mathbb{Q}}
\def\C{\mathbb{C}}

\def\F{\mathbb{F}}
\def\Z{\mathbb{Z}}
\def\D{\mathbb{D}}

\newfont{\EUL}{eufm10 scaled 1000}
\def\J{{\hbox{{\EUL J}}}}
\def\k{{\hbox{{\EUL k}}}}
\def\U{{\hbox{{\EUL U}}}}
\def\u{{\hbox{{\EUL u}}}}
\def\p{{\hbox{{\EUL p}}}}
\def\G{{\hbox{{\EUL G}}}}
\def\m{{\hbox{{\EUL m}}}}

\newcommand\su{{\mathfrak s} {\mathfrak u}}
\renewcommand\sl{{\mathfrak s} {\mathfrak l}}

\newcommand\g{\mathfrak g}
\renewcommand\p{\mathfrak p}
\renewcommand\k{\mathfrak k}

\renewcommand\u{\mathfrak u}

\begin{document}

\begin{center}
{\LARGE \bf Real forms and finite order automorphisms of
\vskip 1cm
affine Kac-Moody algebras - an outline of a 
\vskip 1cm
new approach}
\vskip 2 cm
{\large Ernst Heintze} 
\end{center}

\vskip 2 cm

\section{Introduction}
The classification of real forms and finite order automorphisms of affine Kac-Moody algebras has been achieved by the efforts of many people. In particular the works of F. Levstein \cite{L} and G. Rousseau and his collaborators (\cite{B}, \cite{BR}, \cite{R1}, \cite{R2}, \cite{R3}, \cite{B$_3$R}, \cite{BMR}) have to be mentioned here, but see also {\cite{A}, \cite{Bat}, \cite{BP}, \cite{C}, \cite{JZ}, \cite{Kob} and other papers. The classification probably fills some hundred pages and took about 15 years to get completed.

The purpose of this note is to report on a simpler, quite elementary  approach which in addition gives more complete results. It moreover has the advantage to work in the smooth as well as in the algebraic category, that is for affine Kac-Moody algebras which are extensions of loop algebras consisting of smooth resp. algebraic loops.

While the above mentioned authors always worked in the algebraic setting we are mainly interested in the smooth case which is more appropriate for the purpose of geometry. Actually our interest in these questions orginated from the study of symmetric spaces related to affine Kac-Moody groups and hence from the classification of involutions of ``smooth'' affine Kac-Moody algebras (\cite{HPTT}, \cite{H}). But it turns out that the results are the same in both cases.

Our work started several years ago and, at in early stage, in collaboration with Christian Groß.

This is an expanded version of a talk given at the Symposium ''Geometry related to the theory of integrable systems`` at RIMS, Kyoto, September 2007. Details will appear elsewhere.

\section{Smooth and algebraic affine Kac-Moody algebras}

Instead of working with abstract affine Kac-Moody algebras we directly consider their so called realizations. These are certain two dimensional extensions of (twisted) loop algebras as follows.

Let $\g$ be a simple Lie algebra over the field $\F=\R$ or $\C$ and assume $\g$ in addition to be compact if $\F=\R$. Let $\sigma \in Aut(\g)$ be an arbitrary automorphism, not necessarily of finite order. Then we call
\[
L(\g,\sigma):=\{u:\R\to\g\mid u(t+2\pi)=\sigma u(t),u\in C^\infty\}
\]
a (twisted) loop algebra, $L(\g):=L(\g,id)$ beeing the untwisted loop algebra. $L(\g,\sigma)$ is a Lie algebra w.r.t the pointwise bracket $[u,v]_0(t):=[u(t),v(t)]$. If $\sigma$ has finite order, say $\sigma^l=id$, then the $u\in L(\g,\sigma)$ satisfy $u(t+2\pi l)=u(t)$ and are thus indeed loops. Usually one changes the parameter in this case by the factor $l$, i.e. replaces $u(t)$ by $\tilde u(t):=u(lt)$, and embeds $L(\g,\sigma)$ in this way into $L(\g)$. But this has the slight disadvantage to depend on $l$ (which  not necessarily needs to be the order of $\sigma$ but could be any multiple of it). Moreover such an embedding does not exist if $\sigma$ has infinite order. But we will see later that any $L(\g,\sigma)$ is isomorphic to a twisted loop algebra $L(\g,\tilde\sigma)$ with $\tilde\sigma$ of finite order.

\vskip 0,5 cm
One may weaken the differentiability condition and consider loops of Sobolev class $H^k,k \ge 1$. Everything in the following  works equally well. But this is not so clear for the smallest loop algebra $L_{alg}(\g,\sigma)$, which is usually considered in algebra. This consists of the so called algebraic loops which
are by definition finite Laurent series of the form
\[
u(t)={\sum\limits_{q\in\Q}} u_qe^{iqt}
\] 
with $u_q\in\g$ (resp. $\g_\C$ if $\F=\R$, where $\g_{\C}$ denotes the complexification of $\g$). The periodicity condition $u(t+2\pi)=\sigma u(t)$ requires $u_q$ to lie in the subalgebra $\{x\in\g\mid \sigma^k x=x$ for some $k\in\N\}$ on which $\sigma$ has finite order. In order to ensure surjectivity of the evaluation map $u\mapsto u(t)$ one is hence forced  to assume $\sigma$ to be of finite order in the algebraic case. If $\sigma^l=id$ then periodicity implies that $u(t)$ is actually of the form
\[
u(t)={\sum_{\vert n\vert\le N}}u_ne^{int/l}\ .
\]
Therefore we let
\[
L_{alg}(\g,\sigma)=\{ u\in L(\g,\sigma)\mid u(t)={\sum\limits_{\vert n\vert\le N}}u_n e^{int/l},N\in\N, u_n\in\g_{(\C)}\}.
\]
The definition does not depend on $l$, one only has to assume $\sigma^l=id$. The same remark as above applies here: by changing the parameter by a factor $l$ one might embed $L_{alg}(\g,\sigma)$ into $L_{alg}(\g):=L_{alg}(\g,id)$ and this is usually done. But for our purposes the above definition is more convinient.

So far we have only considered the loop algebras. The affine Kac-Moody algebra is the following $2$-dimensional extension
\[
\hat L(\g,\sigma)=L(\g,c)+\F c+\F d
\]
with
\begin{eqnarray*}
&[u,v]&=[u,v]_0+(u',v)\cdot c \\            
&[d,u]&=u' \\
&[c,x]&=0 
\end{eqnarray*}
for all $u,v,\in L(\g,\sigma)$ and $x\in\hat L(\g,\sigma)$ where $(u,v)={\int\limits^{2\pi}_0}(u(t),v(t))_0dt$. Here $(,)_0$ denotes the Killing form of $\g$ and $u'$ the derivative of $u$.

One easily checks that $\hat L(\g,\sigma)$ is a Lie algebra. The construction could have been done in two steps by introducing  $\tilde L(\g,\sigma):=L(\g,\sigma)+\F c$ first, with brackets as above. This is a one-dimensional central extension of $L(\g,\sigma)$ defined by the cocycle $\omega(u,v):=(u',v)$. $\hat L(\g,\sigma)$ is then a semidirect product of $\tilde L(\g,\sigma)$ with $\F$.

The derived algebra and the center of $\hat L(\g,\sigma)$ are $\tilde L(\g,\sigma)$ and $\F c$, respectively. $L(\g,\sigma)$ is not a subalgebra of $\hat L(\g,\sigma)$ but rather isomorphic to the quotient $\hat L(\g,\sigma)'/\F c$ of the derived algebra by its center.

The extension of $L_{alg}(\g,\sigma)$ to $\hat L_{alg}(\g,\sigma)$ is defined in the same way and the above remarks also apply in this case. In the following we merely consider $\hat L(\g,\sigma)$ and $L(\g,\sigma)$ and come back to the algebraic case only in the last section.

\section{Isomorphisms between affine Kac-Moody algebras}
An important step in our approach is the description of isomorphisms between affine Kac-Moody algebras. They turn out to have a particularly simple form.

Any isomorphism $\hat\varphi:\hat L(\g,\sigma)\to\hat L(\tilde\g,\tilde\sigma)$ induces an isomorphism $\varphi:L(\g,\sigma)\to L(\tilde\g,\tilde\sigma)$ between the loop algebras. Therefore we begin by studying these first. Simple examples of isomorphisms $\varphi:L(\g,\sigma)\to L(\tilde\g,\tilde\sigma)$ are given by
\[
\varphi u(t)=\varphi_t(u(\lambda(t)))
\]
where $\lambda:\R\to\R$ is a diffeomorphism and $t\mapsto \varphi_t:\g\to\tilde\g$ is a smooth curve of isomorphisms. In order that $\varphi u$ (and similarly $\varphi^{-1}u)$  satisfies the periodicity condition $\varphi u(t+2\pi)=\tilde\sigma\varphi u(t)$ for all $t$ we only have to require
\begin{eqnarray*}
(1)& \lambda(t+2\pi)&=\lambda(t)+\epsilon 2\pi\\
(2)& \varphi_{t+2\pi}&=\tilde\sigma\varphi_t\sigma^{-\epsilon}
\end{eqnarray*}
for some $\epsilon\in\{\pm 1\}$. Condition (1) means that $\lambda$ covers a diffeomorphism  $\bar \lambda$ of the circle and $\epsilon=1$ (resp. $-1$) if $\bar\lambda$ and  hence $\lambda$ are orientation preserving (reversing).

We call such isomorphisms {\em standard} and to be of {\em first (second) kind} if $\epsilon=1 (\epsilon=-1)$.

\begin{theorem}\label{3.1}
Any isomorphism $\varphi:L(\g,\sigma)\to L(\tilde\g,\tilde\sigma)$ is standard.
\end{theorem}

The theorem reduces questions about automorphisms of finite order immediately to finite dimensions. It also shows that $\g$ and $\tilde\g$ have to be isomorphic. Therefore we will assume $\tilde\g=\g$ from now on. But $\sigma$ and $\tilde\sigma$ can be different. The periodicity condition (2) gives the only restriction
\[
\tilde\sigma=\varphi_{t+2\pi}\sigma^{\epsilon}\varphi_t
\]
implying that $[\tilde\sigma]$ and $[\sigma]$ are conjugate in $Aut\g/Int\g$. Note that $Aut\g/Int\g$ is isomorphic to the symmetry group of the Dynkin diagram and thus isomorphic to either $1,\Z_2$ or $S_3$(the symmetric group in three letters) and that hence each element is conjugate to its inverse. Moreover  the conjugacy class of $[\sigma]$ is determined by its order, which can be $1,2$ or $3$.

\vskip 0,5 cm
Conversely if $[\sigma]$ and $[\tilde\sigma]$ are conjugate it is easy to find a smooth curve $\varphi_t$ of automorphism satisfying (2). We thus have:

\begin{corollary}\label{3.2}
$L(\g,\sigma)$ and $L(\g,\tilde\sigma)$ are isomorphic if and only if $[\sigma]$ and $[\tilde\sigma]$ are conjugate in $Aut\g/Int\g$. In particular any twisted loop algebra is isomorphic to one with $\sigma$ of finite order.
\end{corollary}

\begin{remark}\label{3.3}\em{
In connection with real forms (section 5) it is interesting to note that Corollary \ref{3.2} also holds  in case $\g$ a real non compact simple Lie algebra (by the same proof). But in this case $Aut\g/Int\g\cong 1,\Z_2,\Z_2\times\Z_2,\D_4$ (the dihedral group) or $S_4$. Hence the order of $[\sigma]$ in $Aut\g/Int\g$ is not enough in this case to distinguish conjugacy classes.}
\end{remark}

The {\em proof of Theorem~\ref{3.1}} consists of several steps. For simplicity let us assume $\sigma=\tilde\sigma=id$ and $\F=\C$. We then can define $\varphi_t$ by $\varphi_t(x)=\varphi(\hat x)(t)$ for all $x\in\g$ where $\hat x$ denotes the constant loop $\hat x(t)\equiv x$. Now, the main point is to prove the existence of a function $\lambda:\R\to\R$ with
\begin{eqnarray*}
(3) &\varphi(f\cdot u)&=(f\circ\lambda)\cdot\varphi(u)
\end{eqnarray*}
for all $u\in L(\g)$ and smooth $2\pi$-periodic $f:\R\to\R$. In fact, if $x_1,\dots,x_n$ is a basis of $\g$ and $u(t)=\Sigma f_i(t)x_i$, we then  get $\varphi(u)(t)=\Sigma f_i(\lambda(t))\varphi_t(x_i)=\varphi_t(u(\lambda(t)))$ as desired. To prove (3) we first show that for any fixed $u,f$ and $t,\ a:=\varphi(fu)(t)$ and $b:=\varphi u(t)$ are linearly dependent. This follows by observing
\[
ad\ a\ ad\ x_1\dots ad \ x_k\ ad\ b=ad\ b\ ad\ x_1\dots ad\ x_k\ ad\ a
\]
for all $x_i\in\tilde\g$ and $k\in\N$ and then applying a classical theorem of Burnside to obtain $ad\ a\ A\ ad\ b=ad\ b\ A\ ad\ a$ for all $A\in End\tilde\g$. We next show $\varphi(fu)(t)=\alpha(f)\cdot\varphi u(t)$ for all $u$ and $f$ but $t$ still fixed for some algebra homomorphism $\alpha$ from the set of $2\pi$-periodic smooth functions to $\C$. In the last step we prove $\alpha (f)=f(t^*)$ for some $t^*\in\R$ and set $\lambda(t):=t^*$.

We finally consider isomorphism $\hat\varphi:\hat L(g,\sigma)\to\hat L(\g,\tilde\sigma)$ between affine Kac-Moody algebras. Since they preserve the center and the derived algebra they are necessarily of the form
\begin{eqnarray*}
&\hat\varphi&=\mu_1c\\
(4)&\hat\varphi d&=\mu_2 d+u_\varphi+\nu_\varphi c\\
&\hat\varphi u&=\varphi u+\alpha(u)\cdot c
\end{eqnarray*}
where $\mu_1,\mu_2,\nu_\varphi\in\F$ are constants, $u_\varphi\in L(\g,\tilde\sigma),\alpha:L(\g,\sigma)\to \F$ is lineare and $\varphi$ is the induced isomorphism between the loop algebras. From Theorem 1 we have $\varphi u(t)=(u(\lambda(t)))$ where $\varphi_t\in Aut\g$ and $\lambda:\R\to\R$ is a diffeomorphism with $\lambda(t+2\pi)=\lambda(t)+\epsilon 2\pi$ and $\epsilon\in\{\pm1\}$. We call $\hat\varphi$ and $\varphi$ to be of the {\em first (second) kind} if $\varphi$ is of the first (second) kind, i.e. if $\epsilon=1$ (resp. $\epsilon=-1)$.

\begin{theorem}\label{3.3}
If $\varphi$ is induced from $\hat\varphi$ then $\lambda$ is linear, i.e. $\varphi u(t)=\varphi_t(u(\epsilon t +t_0))$ for some $\epsilon\in\{\pm 1\}, t_0\in\R$. Conversely, any such isomorphism $\varphi$ is induced by an isomorphisms $\hat\varphi$ between the affine Kac-Moody algebras and this is essentially unique (up to the choice of $\nu_\varphi$ in (4), which can be arbitrary).

More precisely if $\varphi u(t)=\varphi_t(u(\epsilon t+t_0))$ then the $\hat\varphi$ extending $\varphi$ are precisely the ones satisfying $\mu_1=\mu_2=\epsilon, adu_\varphi=-\epsilon\varphi_t'\varphi_t^{-1}$ and $\alpha(u)=-\epsilon(\varphi u,u_\varphi)$ in (4).
\end{theorem}
\begin{corollary}\label{3.4}
There is a bijection between automorphisms of finite order of $\hat L(\g,\sigma)$ and $L(\g,\sigma)$.
\end{corollary}

In fact, if $\hat\varphi$ has finite order then the induced $\varphi$ has finite order. Conversely if $\varphi$ has finite order then there is precisely one $\hat\varphi$ of finite order extending $\varphi$ namely the one with $\nu_\varphi=-{\epsilon\Vert u_\varphi\Vert^2\over 2}$. The reason for this is that $Aut\hat L(\g,\sigma)$ splits as $\{\hat\varphi\in Aut\hat L(\g,\sigma)
\mid\nu_\varphi=-{\epsilon\Vert u_\varphi\Vert^2\over2}\}\times\{\hat\varphi\mid\hat\varphi=id$ on $L(\g,\sigma)+\F c,\hat\varphi d=d+\nu_\varphi c$ for some $\nu_\varphi\in\F\}$ and the second factor contains no elements of finite order.

\section{Automorphisms of finite order}

From the results of the last sections it follows that classifying conjugacy classes of automorphisms of finite order of $\hat L(\g,\sigma)$ is equivalent to classifying conjugacy classes of automorphisms of finite order of $L(\g,\sigma)$ and the aim of this section is to describe such a classification.

Thus let $\varphi:L(\g,\sigma)\to L(\g,\sigma)$ be of finite order. We know that $\varphi$ has the form $\varphi u(t)=\varphi_t(u(\lambda(t)))$ with $\varphi_{t+2\pi}=\sigma\varphi_t\sigma^{-\epsilon}$ and $\lambda(t+2\pi)=\lambda(t)+\epsilon 2\pi$ for some $\epsilon\in\{\pm 1\}$. After a first conjugation we may assume $\lambda(t)=\epsilon t+t_0$ (with $t_0=0$ if $\epsilon=-1$). This comes from the fact that diffeomorphisms of the circle of finite order (like the one induced by $\lambda$) are conjugate to a rotation or a reflection.

Thus we assume $\varphi u(t)=\varphi_tu(\epsilon t+t_0)$. A particularly simple case is the one where $\varphi_t\equiv \varphi_0$ is constant and one may ask whether $\varphi$ is always conjugate to such an automorphism. We have studied those automorphisms in \cite{HPTT}. The answer is the following.

\begin{theorem}\label{4.1}\
\begin{enumerate}
\item [(i)] Not every automorphism of $L(\g,\sigma)$ of finite order is conjugate to one with $\varphi_t$ constant.
\item [(ii)]But for every $\varphi\in Aut(L(\g,\sigma)$ of finite order there exists a $\tilde \sigma\in Aut(\g)$ together with an isomorphismus $\psi:L(\g,\sigma)\to L(\g,\tilde\sigma)$ such that $\tilde\varphi:=\psi\varphi\psi^{-1}$ has constant $\tilde\varphi_t$, that is $\tilde\varphi u(t)=\tilde\varphi_0(u(\epsilon t+t_0))$.
\end{enumerate}
\end{theorem}

We call $\varphi$ and $\tilde\varphi$ {\em quasiconjugate} in the above situation to emphasize that  $\sigma$ has maybe changed. But for automorphisms $\varphi,\tilde\varphi:L(\g,\sigma)\to L(\g,\sigma)$ on the same loop algebra,  quasiconjugate and conjugate are the same.

By the above results (quasi)conjugacy classes of automorphisms of finite order of $L(\g,.)$ and $\hat L(\g,.))$ are classified by certain quadruples $(\epsilon,t_0,\varphi_0,\sigma)$ with $\epsilon\in\{\pm 1\},t_0\in\R$ and $\varphi_0,\sigma\in Aut\g$ modulo some quivalence relation. The equivalence relation of course reflects the fact that certain quadruples correspond to the same quasiconjugacy class.

Although it is possible to prove Theorem \ref{4.1} directly we choose a slightly different path and associate first to any automorophism of finite order an ``invariant'', not only to those with constant $\varphi_t$. We then prove that this is indeed invariant under quasiconjugations and that it moreover distinguishes quasiconjugacy classes. Finally we show that each possible invariant is attained, even by an automorphism with constant  $\varphi_t$. This proves also Theorem~\ref{4.1}, part (ii). Considering automorphisms of order $q$ of the first $(\epsilon=1)$ and second kind $(\epsilon=-1)$ separately, we obtain more precisely the following results.

\begin{theorem}\label{4.2}
Let $\g$ be a compact or complex simple Lie algebra as above. Then the quasiconjugacy classes of automorphisms of the first kind of order $q$ on the various $L(\g,.)$ (or $\hat L(\g,.))$ are in bijection with $\J^q_1(\g):=\{(p,\varrho,[\beta])/p\in\Z,0\le p\le q/2,\varrho\in Aut\g$ from a list of representativs of conjugacy classes of automorphisms of $\g$ of order $(p,q)$ and $\beta\in(Aut\g)^{\varrho}\}$.
\end{theorem}

Here $(p,q)$ is the greatest common divisor of $p$ and $q$ and $[\beta]$ denotes the conjugacy class of $\bar\beta\in\pi_0((Aut\g)^\varrho)=(Aut\g)^\varrho/((Aut\g)^\varrho)_0$.

If $\varphi u(t)=\varphi_t u(t+t_0)$ on $L(\g,\sigma)$ is of order $q$ then we define its invariant as follows. Necessarily $t_0={p\cdot 2\pi\over q}$ for some $p\in\Z$ and we may assume $0\le  p<q$. Let $r:=(p,q),p':=p/r,q':=q/r$ and $l,m\in\Z$ with $lp'+mq'=1$ and $0\le l<q'$. Then $\varphi^{q'}u(t)=\rho_t(u(t))$ and $\varphi^l u(t)=\Lambda_t(u(t+{2\pi\over q'}))$ for some $\rho_t,\Lambda_t\in Aut\g$. Moreover $\rho_t$ has order $r$ and is thus of the form $\rho_t=\alpha_t\varrho\alpha_t^{-1}$ with $\varrho$
from a list of order $q$ automorphisms (modulo conjugation) and $\alpha_t\in Aut\g$. We let $(p,\varrho,[\alpha^{-1}_{t+2\pi/q'}\Lambda^{-1}_t\alpha_t])$ be the invariant of $\varphi$. It is easily checked that this invariant does not change if $\varphi$ is quasiconjugated by an isomorphism of the first kind. But $p$ has to be replaced by $p':=q-p\ (p'=0$ if $p=0$) after a quasiconjugation by an isomorphism of the second kind, which explains the restriction $0\le p\le q/2$ in the definition of $\J_1^q(\g)$. In this way we get a mapping from quasiconjugacy classes to $J^q_1(\g)$. While injectivity of this mapping is the hard part surjectivity follows easily. In fact, the invariant
$(p,\varrho,[\beta])$ may be realized as follows. Let $\sigma:=\varrho^l\beta^{q'},\varphi_0:=\varrho^m\beta^{-p'}$ and $\varphi u(t)=\varphi_0(u(t+p/q 2\pi))$ with $l,m,p',q'$ as above. Then $\varphi$ leaves $L(\g,\sigma)$ invariant, is of order $q$ and has invariant $(\varrho,p,[\beta])$.

To determine conjugacy classes of automorphisms of order $q$ on a fixed $L(\g,\sigma)$ we only have to restrict the invariants to those for which the above example is defined on a loop algebra isomorphic to $L(\g,\sigma)$. By virtue of Corollary \ref{3.2} this is equivalent to $[\varrho^l\beta^{q'}]$ having the same order as $[\sigma]$ in $Aut\g/Int\g$.

In the case of involutions $(q=2)$ things are of course easier. Here $p=0$ or $p=1$. In the latter case $(p,q)=1$, hence $\varrho=id$ and these automorphisms are classified by conjugacy classes of $Aut\g/Int\g$. They are represented by $\varphi u(t)=\varphi_0(u(t+\pi))$ on $L(\g,\varphi_0^{-2})$ where $\varphi_0\in Aut\g$ runs through a list of representations of $Aut\g/Int\g$ (and may thus be chosen to have order $1,2$ or $3$). In particular they do not occur on $L(\g,\sigma)$ if $\sigma^2$ is inner, but $\sigma$ not. In the former case the involutions with invariant $(0,\varrho,[\beta])$ may be represented by $\varphi u(t)=\varrho u(t)$ on $L(\g,\beta)$ where $\varrho\in Aut\g$ runs through the conjugacy classes of involutions and $\beta\in (Aut\g)^\varrho$ runs through the conjugacy class of $(Aut\g)^\varrho/((Aut\g)^\varrho)_0$. This last group is known to be isomorphic to $1,\Z_2,\Z_2\times\Z_2,\D_4$ or $S_4$ where $\D_4$ and $S_4$ denote the dihedral and symmetric groups, respectively. Thus these involutions correspond to finite dimensional symmetric spaces plus a certain extra information.

The situation for automorphisms of the second  kind is described in the following Theorem. Note that their order is necessarily even.

\begin{theorem}\label{5}
Let $\g$ be as above. Then the quasiconjugacy classes of automorphisms of order $q$ of the second kind on the various $L(\g,.)$ are in bijection with
$\J^q_{-1}(\g)=\{(\varphi_+,\varphi_-)\in(Aut\g)^2\mid\varphi^2_+=\varphi^2_-$, order $\varphi_\pm^2=q/2\}/\sim$ where the equivalence relation is generated by $(\varphi_+,\varphi_-)\sim(\varphi_-,\varphi_+)$ and $(\varphi_+,\varphi_-)\sim(\alpha\varphi_+\alpha^{-1},\beta \varphi_-\beta^{-1})$ for all $\alpha,\beta\in Aut\g$ with $\alpha^{-1}\beta\in((Aut\g)^{\varphi_+^2})_0$.
\end{theorem}

If $\varphi\in Aut L(\g,\sigma)$ has order $q$ and is already of the form $\varphi u(t)=\varphi_t(u(-t))$ then $\varphi_t\varphi_{-t}$ has order $q/2$ and there exists $\alpha_t$ with $\varphi_t\varphi_{-t} =\alpha_t\varphi^2_0\alpha_t^{-1}$. The periodicity condition $\varphi_{t+2\pi}=\sigma\varphi_t\sigma$ implies $(\varphi_\pi\sigma^{-1})^2=\varphi_\pi\varphi_{-\pi}$. Hence $\varphi_+:=\alpha_0^{-1}\varphi_0\alpha_0$ and $\varphi_-:=\alpha_\pi^{-1}\varphi_\pi\sigma^{-1}\alpha_\pi$ have order $q/2$ and satisfy $\varphi_+^2=\varphi^2_-$. We thus can define $[\varphi_+,\varphi_-]$ to be the invariant of $\varphi$ and hence get a mapping from the quasiconjugacy classes to $J^q_{-1}(\g)$. Again, surjectivity follows easily. In fact,
the equivalence class $[\varphi_+,\varphi_-]$ may be represented by the automorphismus $\varphi$ on $L(\g,\sigma)$ with $\varphi u(t):=\varphi_+u(-t)$ and $\sigma:=\varphi^{-1}_-\varphi_+$. (Note that $\varphi_+=\sigma\varphi_+\sigma$ and hence $\varphi$ satisfies the periodicity condition).

Let us again look at the special case of involutions $(q=2)$. Since
\[
\J_{-1}^2=\{[\varrho_+,\varrho_-]\mid \varrho_\pm^2=id\}
\]
involutions of the second kind correspond essentially to pairs of compact symmetric spaces. More precisely if $\F=\R$ and $G$ is the compact simply connected Lie group with Lie algebra $\g$ then $[\varrho_+,\varrho_-]$ gives rise to the symmetric pairs $(G,K_+)$ and $(G,K_-)$ where $K_\pm=G^{\varrho_\pm}$ (assuming here $\varrho_{\pm}\ne id)$. The action of $K_+\times K_-$ on $G$ by $(k_+,k_-).g=k_+gk_-^{-1}$ is hyperpolar, that is admits a flat section (actually a torus) which meets every orbit and always orthogonally. The $K_+\times K_-$ action on $G$ 
 is also called a Hermann action and Kollross \cite{K} has shown that essentially all hyperpolar actions on $G$ are Hermann actions. Moreover $(\varrho_+,\varrho_-)$ and $(\tilde\varrho_+,\tilde\varrho_-)$ are equivalent by the above equivalence relation if and only if the corresponding Hermann actions are equivalent. Thus quasiconjugacy classes of involutions of the second kind on the affine Kac-Moody algebras $\hat L(\g,.)$  essentially correspond to equivalence classes of Hermann actions on $G$.

\section{Real forms and involutions}

Let $\g$ be a complex  simple Lie algebra. Real forms of $\hat L(\g,\sigma)$ or $L(\g,\sigma)$ are in bijection with antilinear involutions that is with involutions satisfying $\hat\varphi(ix)=-i\hat\varphi(x)$ (resp. $\varphi(ix)=-i\varphi(x))$. But antilinear automorphismus of finite order on these algebras can be treated exactly the same way as linear ones.

In particular we only have to study the antilinear involutions on $L(\g,\sigma)$ and these
are like in the linear case of the form $\varphi u(t)=\varphi_t(u(\lambda(t)))$, but now with $\varphi_t$ antilinear. Therefore we can associate invariants to them as above and show that their quasiconjugacy classes are parametrized by the sets $\bar\J_\epsilon^q(\g)$ of these invariants where $q$ is the order of $\varphi$ and $\epsilon$ is either $1$ or $-1$ depending on whether $\lambda$ is orientation preserving or reversing. The real forms corresponding to $\bar\J^q_\epsilon(\g)$ are also called {\em almost compact} if $\epsilon=1$ and {\em almost split} if $\epsilon=-1$.
If $\u\subset\g$ is a compact $\sigma$-invariant real form of $\g$ then $\U:=L(\u,\sigma)$ is a real form of $\G:=L(\g,\sigma)$ which is also called a {\em compact} real form. It is the fixed point set of $\varphi u(t)=\omega(u(t)),\omega$ the conjugation of $\g$ w.r.t $\u$.
Automorphisms of order $q$ on $\U$ can be extended to linear as well as to antilinear automorphisms on $\G$ (of order $q$  if $q$ is even and $2q$ if $q$ is odd) and the induced mappings on the sets of invariants turn out to be bijections. For example in case $q=2$ one gets bijections $\J^2_1(\u)\cup \J^1_1(\u)\leftrightarrow  \bar\J_1^2(\g),\J^2_{-1}(\g)\leftrightarrow \bar\J_{-1}^2(\g)$ and $\J_\epsilon^2(\u)\leftrightarrow \J_\epsilon^2(\g)$. To explain these results we recall first the finite dimensional situation.

Up to conjugation the complex, simple Lie algebra $\g$ has precisely one compact real form, say $\u$. If $\varrho$ is an involution of $\u$ then the eigenspace decomposition $\u=\k+\p$ gives rise to the real form $\u^*=\k+i\p$, which is non compact. In this way one gets a bijection between conjugacy class of involutions on $\u$ and conjugacy classes of non compact real forms of $\g$. Moreover each noncompact real form $\g_{\R}$ has by construction a Cartan decomposition, that is a decomposition $\g_{\R}=\k+\m$ into the eigenspaces of an involution such that $\k+i\m$ is a compact Lie algebra and it turns out that this is unique up to conjugation. 

Now the above results may be summerzied as follows.

\begin{theorem}\label{7}

The situation for affine Kac-Moody algebras is exactly the same as that for finite dimensional simple Lie algebras. More precisely:
\begin{enumerate}
\item [(i)] Each complex affine Kac-Moody algebra $\hat\G:=\hat L(\g,\sigma)$ has a ``compact real form'', e.g. $\hat\U:=\hat L(\u,\sigma),\u\subset\g$ a $\sigma$-invariant compact real form, and this is unique up to conjugation.
\item [(ii)]Conjugacy classes of non compact real forms of $\hat\G$ are in bijection with conjugacy classes of involutions of $\hat\U$ (and $\hat\G)$.
\item [(iii)]Each non compact real form of $\hat\G$ has a Cartan decomposition and this is unique up to conjugation.
\end{enumerate}
\end{theorem}

Thus the classification of involutions of $\hat L(\g,\sigma)$ gives also a classification of real forms. 
To describe this more explicitly we first specialize our classification of finite order automorphisms to involutions.

\begin{theorem}\label{5.2}
Up to quasiconjugation the involutions on the various $L(\g,.)$ are given by
\begin{enumerate}
\item [1a)] $u(t)\mapsto \varrho(u(t))$ on $L(\g,\beta)$
\item [1b)] $u(t)\mapsto\varphi(u(t+\pi))$ on $L(\g,\varphi^2)$
\item [2)] $u(t)\mapsto\varrho_+(u(-t))$ on $L(\g,\varrho_-\varrho_+)$,
\end{enumerate}
where $\varrho$ runs through a list of conjugacy classes of involutions of $\g$,  $\beta$ and $\varphi$ represent conjugacy classes of $\pi_0((Aut\g)^\varrho),\pi_0(Aut\g)$ resp. and $(\varrho_+,\varrho_-)$ represent equivalence classes of $\{(\varrho_+,\varrho_-)\in(Aut\g)^2/\varrho_{\pm}^2=id\}/\sim$. Since all automorphisms can be chosen to leave the compact real form $\u$ of $\g$ invariant, the above mappings classify also involutions on the various $L(\u,.)$.
\end{theorem}

The invariants of these involutions are $(0,\varrho,[\beta]),(1,id,[\varphi])$ and $[\varrho_+,\varrho_-]$, respectively. The involutions of type 1b) correspond to the various isomorphism classes of the $L(\g,\sigma)$. If $\g$ has no  outer autmorphism there exists
precisely one such involution (given by $u(t)\mapsto u(t+\pi)$ on $L(\g))$ and otherwise two unless $\g$ is of type $D_4$ which has threee.

According to Theorem \ref{7} the real forms (of the various  $L(\g,.)$ and up to quasiconjugation) are thus the {\em compact} real forms $L(\u,.)$ (corresponding to $J_1^1(\u))$ and the {\em fixed point sets} of $\tilde\omega\circ\psi$ where $\psi$ is an involution from the list above. Here $\tilde\omega(u)(t)=\omega(u(t)),$ where $\omega$ denotes conjugation w.r.t. $\u$. From this we obtain the following  classification.

\begin{theorem}\label{5.3}
The real forms of the  various $L(\g,.)$ are up to quasiconjugation represented by
\begin{enumerate}
\item [1a)] $L(\g_\R,\sigma)\subset L(\g,\sigma)$ where $\g_\R$ runs through all conjugacy classes of real forms of $\g$ and $\sigma\in Aut\g_\R$ runs through a list of representatives of conjugacy classes of $\pi_0(Aut\g_\R)$.
\item [1b)] $\{u\in L(\u,\varphi^2)\mid u(t+\pi)=\varphi(u(t))\}+i\{u\in L(\u,\varphi^2)\mid u(t+\pi)=-\varphi(u(t))\}$, where $\varphi\in Aut\u$ represents a conjugacy class of $\pi_0(Aut\u)\cong\pi_0(Aut\g)$.
\item[2)] $\{\Sigma u_ne^{int/l}\in L(\g,\sigma)\mid u_n\in\g^{\omega\varrho_+}\}$ where $\sigma=\varrho_-\varrho_+$ has order $l$ and $(\varrho_+,\varrho_-)$ represents an equivalence class of $\{(\varrho_+,\varrho_-)\in(Aut\g)^2\mid\varrho_+^2=\varrho^2_-\}/\sim$.
\end{enumerate}
The corresponding real forms of $\hat L(\g,.)$ are obtained by adjoining $\R c+\R d$ in case 1a) and 1b) and $\R(ic)+\R(id)$ in case 2).
\end{theorem}
\begin{remarks}\em{
\hspace{3pt}
\begin{enumerate}
\item[(i)] $\varrho_\pm$ can always be chosen such that $\varrho_-\varrho_+$ has finite order (in fact $\le 4$), cf. also Lemma \ref{6.2}.
\item[(ii)] While the real forms in 1a) and 1b) constitute the almost compact ones, those in 2) are the almost split.
\item[(iii)] Cartan decompositons can be obtained in all cases above by intersecting the real forms with curves in $\u$ and $i\u$, respectively (assuming that $\g_\R$ is $\omega$-invariant in 1a)). For example, if one considers in 2) only loops $\Sigma u_n\lambda^n:=\Sigma u_n e^{int/l}$ which converge in the whole puncutred plane $\C^*$ then the two summands of the Cartan decomposition consist of functions whose restriction to the real line is contained in $\g_\R:=\g^{\omega\varrho_+}$ while their restriction to the unit circle is contained in $\u$ and $i\u$, respectively (and which satisfy the appropriate periodicity condition).
\end{enumerate}}
\end{remarks}

To illustrate Theorem {\ref {5.3}} we consider the simplest case $\g=\sl(2,\C)$. The algebra $\sl(2,\C)$ has no outer automorphisms and up to conjugation only one involution, which may be represented e.g. by 
$\tau=Ad\left(
\begin{array}{cc}
1&\\
&-1
\end{array}\right)$
or  by $\mu$ with $\mu A=-A^t$. As compact real form we  may take $\u=\su(2)$. We let $\omega$ be the conjugation w.r.t. $\u$, i.e. $\omega A=-\bar A^t$. In particular, $\g^{\omega\mu}=\sl(2,\R)$. Then the real forms of $L(\sl(2,\C))$ are up to quasiconjugation the following.
\begin{enumerate}
\item[1a)] $L(\su(2))$, $L(\sl(2,\R))$, $L(\sl(2,\R),\tau)$
\item[1b)] $\{u\in L(\u)\mid u(t+\pi)=u(t)\}+i\{u\in L(\u)\mid u(t+\pi)=-u(t)\}$
\item[2)] $\{\Sigma u_ne^{int}\in L(\sl(2,\C))\mid u_n\in\su(2)\}$, $\{\Sigma u_ne^{int}\in L(\sl(2,\C))\mid u_n\in\sl(2,\R)\}$, and $\{\Sigma u_ne^{int/2}\in L(\sl(2,\C),\mu)\mid u_n\in\sl(2,\R)\}$.
\end{enumerate}
They correspond to the invariants $(0,id,[id])\in \J^1_1(\sl(2,\C)),(0,\mu,[id]),(0,\mu,[\tau]),(1,id,[id])$
$\in \J^2_1 (\sl(2,\C))$ and $[id,id],[\mu,\mu]$, $[\mu,id] \in\J_2^2(\sl(2,\C))$, respectively. They are all real forms of $L(\sl(2,\C))$ except the third and the last one which are real forms of  a twisted loop algebra of $\sl(2,\C)$. But an isomorphism $\psi:L(\sl(2,\C),\sigma)\to L(\sl(2,\C))$ is for example given by $\psi u(t)=\psi_t(u(t))$ with $\psi_t=e^{ad tX}$ such that $\psi_{2\pi}=\sigma^{-1}$ and this carries the corresponding real form into one of $L(\sl(2,\C)$.

\section{The algebraic case} \label{6}

In contrast to the above smooth setting all authors so far (with the exception of \cite{HPTT}) have  considered the algebraic case, i.e. automorphisms and real forms of $\hat L_{alg}(\g,\sigma)$ and $L_{alg}(\g,\sigma)$. This case is more rigid and thus more subtle. For example it is much harder to find and ``algebraic'' isomorphism which conjugates two given  smoothly conjugate automorphisms of $\hat L_{alg}(\g,\sigma)$. But it turns out that our methods and ideas also work in this setting when  suitably refined and combined with a  result of Levstein \cite{L}. In fact we have:

\begin{theorem} \label{8}

Conjugacy classes of finite order automorphisms and real forms of $\hat L_{alg}(\g,\sigma)$ and $L_{alg}(\g,\sigma)$ are classified by the same invariants as in the smooth case and are thus in bijection to those of $\hat L(\g,\sigma)$ or $L(\g,\sigma)$, respectively. Moreover the noncompact real forms of $\hat L_{alg}(\g,\sigma)$ admit a Cartan decomposition and this is unique up to conjugation.
\end{theorem}

The proof follows along the same lines as above but has to be modified at several points. 

To begin with, isomorphisms $\varphi:L_{alg}(\g,\sigma)\to L_{alg}(\g,\tilde\sigma)$ are not necessarily standard (assuming $\F=\C$ in the following). They are rather of the form $\tilde\varphi\circ\tau_r$ with $\tilde\varphi$ standard and
\[
\tau_r({\sum\limits_{\vert n\vert\le N}u_ne^{int/l}})={\sum\limits_{\vert n\vert\le N}u_nr^n e^{int/l}}
\]
for some positive $r$. The reason for this is that the algebra homomorphism $\alpha:C^\infty(S^1)\to\C$ in the proof of Theorem \ref{3.1} is not necessarily evaluation at some point when restricted to the algebraic functions $\{{\sum\limits_{\vert n\vert\le N}a_ne^{int}}\}$ but more generally of the form ${\sum\limits_{\vert n\vert\le N}}a_nz^n\to{\sum\limits_{\vert n\vert\in N}}a_nz_0^n={\sum\limits_{\vert n\vert\le n}}a_nr^ne^{{int}^*}$. Fortunately, it turns out that these extra isomorphisms $\tau_r$ do not affect the discussion of finite order automorphisms in an essential way. For example automorphism of the first kind of finite order are always standard and those of the second kind are conjugate to a standard one. Thus it is completely sound to work only with standard isomorphisms. Since these can be viewed as special isomorphisms between the corresponding smooth algebras we can apply the results of section 4 and associate to each autmorphism of order $q$ of $L_{alg}(\g,\sigma)$ an invariant and thus get a mapping $I^q_\epsilon$ from the set $Aut_\epsilon^q(L_{alg}(\g,.))/Aut(L_{alg}(\g,.))$ of quasiconugacy classes of such automorphisms on the various $L_{alg}(\g,.)$ into $\J^q_\epsilon(\g)$, where $\epsilon=1$ or $-1$ depending on whether the automorphisms are of the first or second kind. The goal is to prove bijectivity of these mappings.

Surjectivity follows essentially from what has been done in the smooth case. In fact, the construction there yielded for each invariant a $\sigma\in Aut\g$ and an automorphism $\varphi$ of $L(\g,\sigma)$ with this invariant of the form $u(t)\mapsto\varphi_0(u(\epsilon t+t_0))$. This $\varphi$ also preserves $L_{alg}(\g,\sigma)$, but in the algebraic case we have to ensure that $\sigma$ has finite order. For example if $[\varphi_+,\varphi_-]\in \J^q_{-1}(\g)$ then the construction yielded $\sigma:=\varphi_-^{-1}\varphi_+$ which is in general not of finite order. But $(\varphi_+,\varphi_-)$ can replaced by an equivalent pair and the result follows by showing first the existence of a $\varphi_\pm$-invariant compact real form $\u$ of $\g$ and then by applying the next Lemma to $G:=\{\alpha\in(Aut\g)^{\varphi_\pm^2}\mid\alpha(\u)=\u\}$.

\begin{lemma} \label{6.2} Let $G$ be a compact Lie group and $g_+,g_-\in G$ with $g_+^2=g_-^2$. Then there exists $h\in G_0$ such that $(hg_-h^{-1})^{-1}\cdot g_+$ is of finite order.
\end{lemma}

But the main problem is the injectivity of $I_\epsilon^q$. Here our elementary methods do not suffice. To solve it, we use the following basic result of Levstein.

\begin{theorem} {\rm (Levstein \cite{L}, \cite{R2})} Let $\hat\varphi$ be an automorphism of $\hat L_{alg}(\g,\sigma)$ of finite order. Then $\hat L_{alg}(\g,\sigma)$ has a Cartan subalgebra which is invariant under $\hat\varphi$.
\end{theorem}

Since after a conjugation we may assume that $\hat\varphi$ leaves the standard Cartan subalgebra invariant which consists of constant loops $+\C c+\C d$ we conclude that the $u_\varphi$ from $\hat\varphi d=\epsilon d+u_\varphi+\nu_\varphi c$ is constant and this implies $\varphi_t=e^{adtX}\varphi_0$ for some $X\in\g$ and $\varphi_0\in{Aut}\g$, because $\varphi_t'\varphi_t^{-1}=-\epsilon\ adu_\varphi$ (assuming $\varphi u(t)=\varphi_t(\epsilon t+t_0))$. In the next step we even get rid of the $e^{adtX}$-factor by a further quasiconjugation. Here it is essential to allow $\sigma$ to be replaced by some $\tilde\sigma\in Aut\g$. Thus we have also in the algebraic case the result that any automorphism of finite order is quasiconjugate to one of the form $\hat\varphi c=\epsilon c,\hat\varphi d=\epsilon d$, and $\varphi u(t)=\varphi_0(u(\epsilon t+t_0))$. The proof of the injectivity of $I_\epsilon^q$ is therefore considerably simplified in that we have only to consider  these very special automorphisms. Surprisingly, it follows finally from the hyperpolarity (s. section 4) of the $\sigma$-action (for automorphisms of the first kind) and the Hermann action (for automorphisms of the second kind). If $G$ is a compact connected Lie group and $\sigma\in Aut\ G$ then the action of $G$ on itself by $g.x:=gx\sigma(g)^{-1}$ is called the $\sigma$ action. A maximal torus of $G^\sigma=\{g\in G/\sigma g=g\}$ meets every orbit and always orthogonally (w.r.t. any biinvariant metric).

\vskip 0,5 cm
The same ideas also work for antilinear automorphisms of finite order and in the real case $(\F=\R)$. Therefore the results about real forms and Cartan decompositions carry over from the smooth to the algebraic setting without any difficulties.

\eject

\vskip3cm
Ernst Heintze\\
Institut für Mathematik, Universität Augsburg\\
Universitätsstrasse 14\\
D - 86159 Augsburg, Germany\\
e-mail address: {\it heintze@math.uni-augsburg.de}\\

\end{document}